\documentclass[11pt]{article}
\usepackage{amssymb,amsmath,amsfonts,amsthm,mathtools,color}

\usepackage{mathrsfs}
\usepackage{bm} % for bm macro

\usepackage{comment} % for comment large section
\usepackage[title,titletoc,header]{appendix}
\usepackage{graphicx,subfig}

\usepackage{paralist}
\usepackage{pgfplots}
\usepackage{tikz}
\usetikzlibrary{arrows,positioning,shapes.geometric}

\usepackage[ruled,vlined]{algorithm2e}

\graphicspath{ {./figures/} }

\usepackage[utf8]{inputenc}

\usepackage{indentfirst}
\usepackage{multicol}
\usepackage{booktabs}
\usepackage{url}
\usepackage[outdir=./]{epstopdf} %to enable the use of eps
\usepackage{stmaryrd}
\usepackage[shortlabels]{enumitem}
\usepackage{hyperref}
\hypersetup{
    colorlinks=true, %set true if you want colored links
    linktoc=all,     %set to all if you want both sections and subsections linked
    linkcolor=blue,  %choose some color if you want links to stand out
}
\numberwithin{equation}{section}

\theoremstyle{definition}

\theoremstyle{remark}
\newtheorem{Remark}{Remark}[section]

\def\to{\rightarrow}

\def\del{\partial}

 \def\bs{\bigskip}

\def\bs{\boldsymbol}

\def\cD{\mathcal{D}}

\def\cG{\mathcal{G}}

\def\cL{\mathcal{L}}

\def\cQ{\mathcal{Q}}

\def\cU{\mathcal{U}}

\def\d{{\mathrm{d}}}

\def\sE{{\mathbb{E}}}

\def\sR{{\mathbb R}}
\def\sS{{\mathbb{S}}}

\DeclareMathOperator*{\argmax}{arg\,max}

\newcommand{\lc}
{\mathrel{\raise2pt\hbox{${\mathop<\limits_{\raise1pt\hbox
{\mbox{$\sim$}}}}$}}}

\newcommand{\gc}
{\mathrel{\raise2pt\hbox{${\mathop>\limits_{\raise1pt\hbox{\mbox{$\sim$}}}}$}}}

\newcommand{\ec}
{\mathrel{\raise2pt\hbox{${\mathop=\limits_{\raise1pt\hbox{\mbox{$\sim$}}}}$}}}

\def\bb{\begin{equation}} \def\ee{\end{equation}}
\def\bbn{\begin{equation*}} \def\een{\end{equation*}}

\def\beqn{\begin{eqnarray}}  \def\eqn{\end{eqnarray}}

\def\beqnx{\begin{eqnarray*}} \def\eqnx{\end{eqnarray*}}

\def\bn{\begin{enumerate}} \def\en{\end{enumerate}}

\def\bd{\begin{description}} \def\ed{\end{description}}



\begin{document}

\title{
Generalised correlated batched bandits via the ARC algorithm with application to dynamic pricing}

\author{
	Samuel N. Cohen\footnote{Mathematical Institute, University of Oxford and Alan Turing Institute,   \texttt{samuel.cohen@maths.ox.ac.uk}}
	\and
	Tanut Treetanthiploet\footnote{Alan Turing Institute,   \texttt{ttreetanthiploet@turing.ac.uk}}
}
\date{\today}

\maketitle
\begin{abstract}
The Asymptotic Randomised Control (ARC) algorithm provides a rigorous approximation to the optimal strategy for a wide class of Bayesian bandits, while retaining low computational complexity. In particular, the ARC approach provides nearly optimal choices even when the payoffs are correlated or more than the reward is observed. The algorithm is guaranteed to asymptotically optimise the expected discounted payoff, with error depending on the initial uncertainty of the bandit. In this paper, we extend the ARC framework to consider a batched bandit problem where observations arrive from a generalised linear model. In particular, we develop a large sample approximation to allow correlated and generally distributed observation. We apply this to a classic dynamic pricing problem based on a Bayesian hierarchical model and demonstrate that the ARC algorithm outperforms alternative approaches.

		\smallskip
		
		
\end{abstract}
\section{Introduction}

In the multi-armed bandit problem, a decision maker needs to sequentially decide between acting to reveal data about a system and acting to generate
profit. The central idea of the multi-armed bandit is that the agent has $K$ `options' or `arms', and must choose which arm to play at each time. Playing an arm results in a reward generated from a fixed but unknown distribution which must be inferred on-the-fly. 

In the classic multi-armed bandit problem, the reward of each arm is assumed to be independent of the others  (see \cite{Gittin_origin, Original_UCB, Lattimore_book})  and it is the only observation obtained by the agent at each step. In practice, we often observe signals in addition to the rewards and there is often correlation between the distributions of outcomes for different choices (i.e.~playing one arm may give information about the other arms).

The differences between signals (observations) and rewards are important for designing learning algorithms, especially when the arms are correlated. In particular, it is possible that a single arm is informative but costly. Most algorithms for bandits will fail to choose this arm, since it never yields the best immediate reward (see \cite{ARC} for discussion).

In \cite{ARC}, we develop the Asymptotic Randomised Control (ARC) approach to the learning problem with general observations by considering the dynamics of the posterior estimates together with a smooth asymptotic approximation. The ARC approach yields a natural interpretation of the value of each arm as a sum between the exploitation gain and the learning premium, as in the UCB principle \cite{UCB_tuned, Auer_UCB_Bandit}. One of the main limitations of the original ARC approach is that it requires the model to have an explicit posterior update, which fundamentally requires the observations and prior to be a conjugate pair. This paper aims to extend this framework to study a more general class of learning problems, where the observations arrive from a Generalised Linear Model (GLM). 

Some works have been developed to study bandits with correlation (see \cite{GLM_bandit, Linear_Gaussian_bandit, Tutorial_on_Thompson_Sampling, Knowledge_Gradient}) and theoretical guarantees are often provided through regret analysis for specific examples. In particular, the reward is assumed to be the only observation and the structural correlation is assumed on the reward.  Filippi et al. \cite{GLM_bandit} considers the correlated bandit problem with GLM rewards (i.e.~the reward of the $i$th arm has mean $\mu(\theta^\top x_i)$, where $\theta$ is an unknown parameter, $x_i$ is a known feature of the $i$th arm and $\mu$ is a known invertible link function. In contrast, we will model generic observations through a GLM, but allow increased flexibility in modelling rewards.

\subsection{Dynamic Pricing Example}
\label{sec:dynamic pricing intro}
To illustrate an application where modelling observations (rather than rewards) is useful, we consider a simple dynamic pricing problem where our agent needs to learn customer demand while maximising revenue. We will refer to this setting later in this paper, to provide a concrete example of our problem setup and our solution. In fact, one can easily extend this framework to study a general class of learning problems. 
% We will refer to this setting later in this paper, to provide a concrete example of our problem setup and our solution via a simulation. 

Suppose that at the beginning of each day, our agent needs to choose a price from the set $\{c_1, ..., c_K\}$. At the end of day, with chosen price $c_k$, our agent observes $N_k$ customers arriving at the store and a sequence of random variables $\big(Q^{(k)}_{i}\big)_{i=1,...,N}$, taking values in $\{0,1\}$, indicating whether the product is bought by the $i$th customer. Assuming independence between the number of customers arriving and their purchasing decisions, the expected revenue of the agent on that day can be given by
$$r(c_k) := \sE\bigg[c_k \sum_{i=1}^{N}Q^{(k)}_{i}\bigg] = c_k \sE[N] p(c_k),$$
where $p(c_k)$ is the probability that each customer will buy the product, given that the price is $c_k$.

In September 2015,  Dub\'e and Misra \cite{Data_dynamic_pricing} ran an experiment, in collaboration with the business-to-business company ZipRecruiter.com, to choose an optimal price for an online subscription. Their experiment ran in two stages, as an offline learning problem: first collecting data using randomly assigned prices, and then testing their optimal price. 

\begin{figure}[h]
	\centering
	\includegraphics[trim=8cm 1cm 8.5cm 1cm, width=0.7\linewidth]{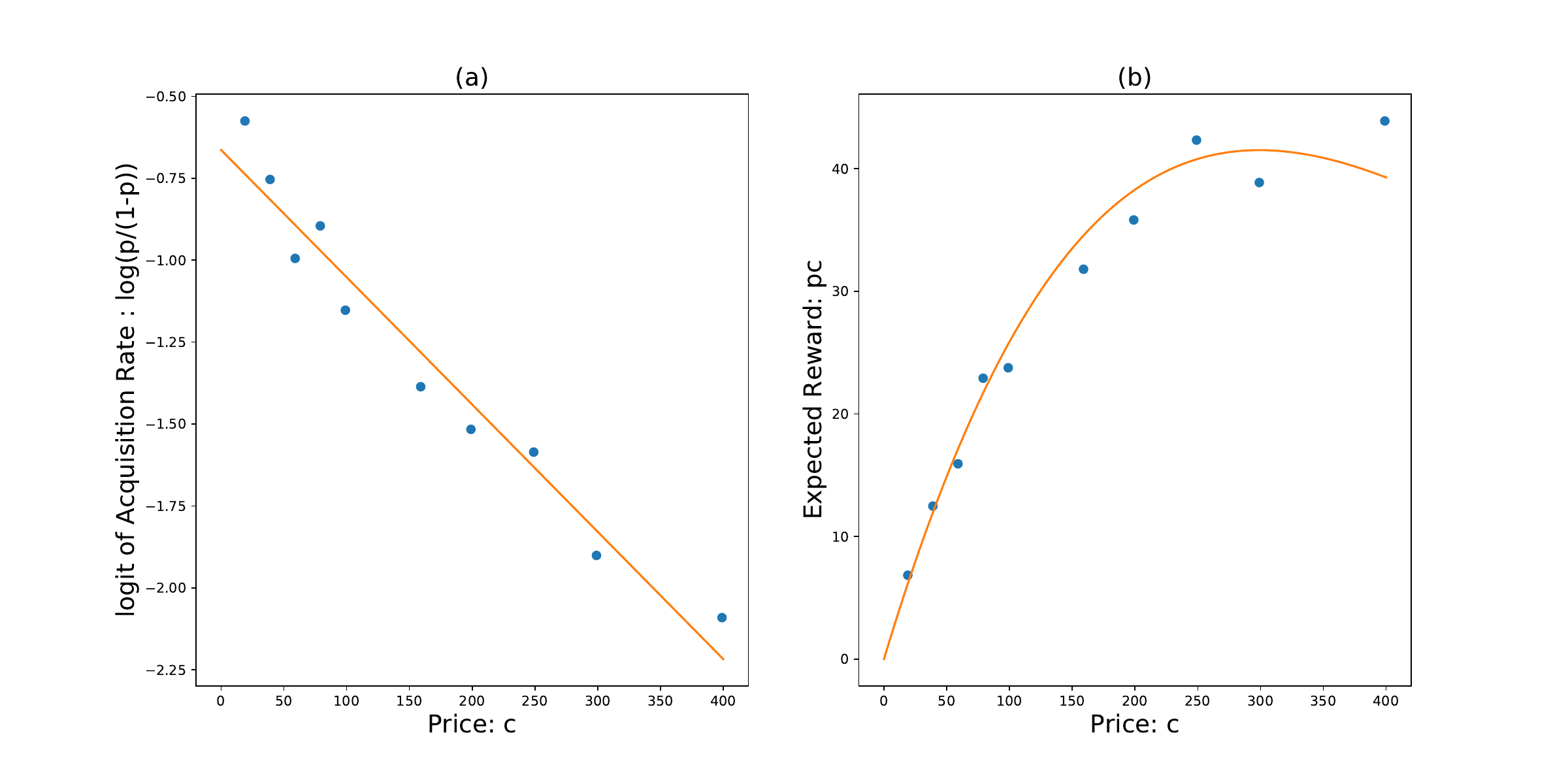}
	\caption{(a): The relation between the price and the logit of Acquisition Rate. (b): The relation between the price and the expected reward per customer, with the implied relation from the regression line in acquisition rate. }
	\label{fig: price logit relation}
\end{figure}

Figure \ref{fig: price logit relation}(a) displays the relation between the logit of the acquisition rate (the proportion of the customers who subscribe), together with its best fit line. It is clear that the demands at different prices are related (i.e.~demand decays with price), which agrees with economic intuition. These data can be found in Table 5 of \cite{Data_dynamic_pricing}. We also observe that an approximately linear relation occurs between prices the logit function of the conversion rate, as seen in Figure \ref{fig: price logit relation}(a), but not the expected reward as observed in Figure \ref{fig: price logit relation}(b). 

Guided by to Figure  \ref{fig: price logit relation}, it is reasonable to consider a logistic model (a case of the Generalised Linear Model (GLM)) for the probability of subscription; the reward, however, does not fit as naturally into a GLM framework. This means that this problem does not easily fit in the correlated bandit frameworks where a distributional assumption is directly made on the reward. 

\subsection{Contribution of this paper}
The novelty of this work is to illustrate a method to study the correlated bandit problem, where observations and rewards are different, and to impose structural assumptions on the observations, rather than the rewards.

In this paper, we will extend the application of the ARC algorithm \cite{ARC} to a wider class of learning problems, where the observations of each arm arrive from a parametric exponential family in batches, with shared parameter, following a GLM. We generalise the implementation of the ARC algorithm to tackle cases without exact prior-posterior conjugacy, by applying the Kalman filter with a large sample observation. We also provide a numerical simulation to illustrate the performance of our modification. 

This paper proceeds as follows. In Section \ref{sec:ARC review}, we give an overview of the ARC approach proposed in \cite{ARC}. We then describe the batched bandit problem, in Section \ref{sec:batched setup}, formulated through the generalised linear bandit model and describe how we can use large sample theory and Bayesian statistics to propagate our posterior. Finally, in Section \ref{sec:Num sim}, we describe how we can use offline data to construct a bandit environment, and use experimental data from \cite{Data_dynamic_pricing} to simulate the pricing problem and compare our method with other approaches.

\section{Frameworks for the ARC approach}
\label{sec:ARC review}
 To avoid unnecessary technical difficulty, we shall describe the setup of \cite{ARC} without describing the full set of technical assumptions under which asymptotic optimality can be proven. The reader should consult the original paper for formal discussion. 

\subsection{Problem setup}
Suppose that our bandit has $K$ arms.  Let $\bs{\theta} \sim \pi$ be an underlying (unknown) parameter taking values in $\Theta$ with a (known) prior $\pi$. The parameter $\bs{\theta}$ describes the distribution of our bandit, i.e.~when the $i$th arm is chosen at time $t$, we observe a random variable $Y^{(i)}_t \sim \pi_i(\cdot | \bs{\theta})$ and obtain a reward $r_i(Y^{(i)}_t)$. 

The prior (and posterior) of $\theta$ is assumed to be described by a finite dimensional process $\zeta_t = (M_t, D_t)$ taking values in a set $\Theta \times \cD \subseteq \sR^p \times \sR^q$. In particular, the posterior of $\bs{\theta}$ at time $t$ is given by $\cL(\bs{\theta} | \cG^U_t)  = \cL(\bs{\theta}|M_t, D_t)$. Here, $\cL(\bs{\theta} | \cG^U_t)$ denotes the law of $\bs{\theta}$ conditional on the historical observations $\cG^U_t$. The process $(M_t)$ represents an estimate (i.e.~posterior mean) of the unknown parameter $\bs{\theta}$, while $(D_t)$ represents statistical error (i.e.~posterior variance) of our estimate. This parameterisation is shown to be useful in the analysis given in \cite{ARC}.

The transition of the posterior/prior when the $i$th arm is chosen is assumed to be updated by
\begin{equation}
\label{eq: transition map}
\zeta_t  := 
\zeta_{t-1} +
\begin{pmatrix}
\mu_i(\zeta_{t-1}) 
 \\ b_i(\zeta_{t-1}) 
\end{pmatrix} +
\begin{pmatrix}
\sigma_i(\zeta_{t-1}) Z^{(i)}_t\\0
\end{pmatrix} ,
\end{equation}
where $\mu_i :  \Theta \times \cD \to \sR^p$,  $b_i:  \Theta \times \cD  \to \cD $, $\sigma_i :  \Theta \times \sR^{q}  \to \sR^{p \times r}$ and $Z^{(i)}_t$ is a random variable obtained from rescaling $Y^{(i)}_t$, such that $\sE[Z^{(i)}_t |M_{t-1}, D_{t-1}] = 0$ and $\text{Var}[Z^{(i)}_t |M_{t-1}, D_{t-1}] = 1$.  The random variable $Z^{(i)}_t$ in the above update appears through the predictive distribution of the new observation.

\subsection{Example: Linear Gaussian bandit}
Suppose that our multi-armed bandit problem associates to the unknown parameter $\bs{\theta}$ where $\bs{\theta} \sim N(m,  d)$ describes the prior of $\bs{\theta}$ where $m \in \sR^p$ and $d \in \sS_+^p$. 

Suppose that when the $i$th arm is chosen, our observation has distribution $Y^{(i)} \sim N(x_i^\top \bs{\theta},  1)$ where $x_i$ is a known vector in $\sR^p$. Following classical Bayesian analysis, we can describe the posterior at time $t$ of $\bs{\theta}$ by $\cL(\bs{\theta} | \cG^U_t) = N(M_t, D_t)$ where $(M_t, D_t)$  follows 
\begin{equation}
\label{eq:posterior dynamic}
\begin{aligned}
M_{t} &= \big(D_{t-1}^{-1} + x_{I_{t}} x_{I_{t}}^\top\big)^{-1}\big(D_{t-1}^{-1} M_{t-1} + x_{I_{t}} Y^{(i)}_{t}\big), \\
D_{t} &= \big(D_{t-1}^{-1} + x_{I_{t}} x_{I_{t}}^\top\big)^{-1},
\end{aligned}
\end{equation}
for $(I_t)$ a process taking values in $\{1, 2, ..., K \}$ indicating the arm we choose at each time.

Since {\small $\cL(Y^{(i)}_{t} | M_{t-1}, D_{t-1}) = N\big(c_i^\top M_{t-1},  c_i^\top D_{t-1} c_i + 1 \big)$}, we can write \eqref{eq:posterior dynamic} in the form of \eqref{eq: transition map}, for $I_t = i$, as
{\small
\begin{equation}
\label{eq: Gaussian mean propagation}
\begin{aligned}
M_{t} & = M_{t-1} +  \big(D_{t-1}^{-1} + x_{i} x_{i}^\top\big)^{-1} \big( x_i^\top D_{t-1} x_i + 1 \big)^{1/2} Z^{(i)}_t, \\
D_t&= D_{t-1} - D_{t-1} x_{i}(1 + x_{i}^\top D_{t-1} x_{i})^{-1}x_{i}^\top D_{t-1},
\end{aligned}
\end{equation}
}
where $Z^{(i)}_{t} := \Big( \frac{Y^{(i)}_t - x_i^\top M_{t-1}}{\sqrt{x_i^\top D_t x_i + 1}} \Big) \sim N(0, 1)$ conditional on $(M_{t-1}, D_{t-1})$.

\subsection{Problem objective and ARC solution}
Assume the agent chooses the arm $I_t$ at time $t$ by sampling from the probability $U_t$ taking values in $\Delta^K := \big\{ u \in [0,1]^K : \sum_{i=1}^K u_i = 1 \big\}$. Using the dynamics \eqref{eq: transition map}, \cite{ARC} study the discounted reward over infinite horizon and rewrite the objective function in terms of a Markov decision process with underlying state $(M_t,D_t)$ and objective $V(m,d) :=
 \sup_{U \in \cU}V^U(m,d)$, where
\begin{equation}
\label{eq:objective}
\begin{aligned}
    &V^U(m,d)  = \sE_{m,d} \bigg[\sum_{t=1}^\infty \beta^{t-1} \Big(\sum_{i=1}^K U_{i,t} f_i(M_{t-1}, D_{t-1}) \Big)\bigg]
    \end{aligned}
\end{equation}
with $f_i(m,d) := \sE_{m,d} \big[r_i(Y^{(i)}) \big]$.

It was shown in \cite{ARC} that an asymptotic solution (for small $d$) of the above problem is to sample using a distribution given by a feedback function (of the posterior parameter) $U_\rho(m,d) = \nu \big({\rho \| d\|}, \alpha({\rho \|d\|}, m,d)\big) \in \Delta^K$,  where:

$\bullet$ $\rho \in (0,\infty)$ is a chosen hyper-parameter. 

$\bullet$ $\nu : (0,\infty) \times \sR^K \to \Delta^K$ is a function given by 
    	\begin{equation}
	    \label{eq:nu}
	    \nu_i(\lambda, a) = \frac{\exp\big(a_i/\lambda\big)}{\sum_{j=1}^K\exp\big(a_j/\lambda\big)}.
	\end{equation}
$\bullet$ $\alpha : (0,\infty) \times \Theta \times \cD \to \sR^K$ is a function given by $\alpha(\lambda, m,d) = f(m,d) + \beta(1-\beta)^{-1}L(\lambda,  m,d)$ with $f:\Theta \times \cD \to \sR^K$  and $L : (0,\infty) \times \Theta \times \cD \to \sR^K$.

The function $\alpha$ represents the (incremental) values of each arm obtained by the summation of the function $f$ describing the exploitation gain and the function $L$ describing the learning premium, which can be given explicitly in terms of the derivatives of $f$ by 
{
\begin{equation}
\label{eq: L expression}
\begin{aligned}
L_i(\lambda, m,d) &:= \langle \mathcal{B}^\lambda(m,d) ; b_i(m, d) \rangle + \langle \mathcal{M}^\lambda(m,d) ; \mu_i(m,d) \rangle +  \frac{1}{2}\text{Tr} \big( \Sigma^\lambda(m,d) \sigma_i\sigma_i^\top(m,d) \big)
\end{aligned}
\end{equation}
}
with $b_i, \mu_i$ and $\sigma_i$ given in \eqref{eq: transition map} and $\mathcal{B}^{\lambda} : \Theta \times \cD  \to \sR^p$,  $\mathcal{M}^{\lambda} : \Theta \times \cD  \to \cD $ and $\Sigma^{\lambda} : \Theta \times \cD  \to \sR^{p \times p}$ satisfying
{\small
\begin{equation}
\label{eq: B, M, Sigma Expression}
\begin{aligned}
	&\mathcal{B}^{\lambda}(m, d) := \sum_{j =1}^K \nu^{\lambda}_j\big( \lambda, f(m,d) \big) \del_d f_j(m,d),  \qquad \mathcal{M}^{\lambda}(m, d) := \sum_{j =1}^K \nu^{\lambda}_j\big(\lambda, f(m,d) \big) \del_m f_j(m,d),  \\
	&\Sigma^{\lambda}(m,d) :=\sum_{j =1}^K \Big(\nu_j\big(\lambda, f(m,d) \big)\del^2_m f_j(m,d)\Big) + \frac{1}{\lambda}\sum_{i,j =1}^K\Big(\eta^{\lambda}_{ij}\big((f(m,d) \big) \big(\del_m f_i(m,d) \big) \big(\del_m f_j(m,d) \big)^\top\Big). 
	\end{aligned}
\end{equation}
}
with $\eta_{ij}^\lambda(a) :=  \nu_i(\lambda, a)(\mathbb{I}(i=j) - \nu_j(\lambda, a) )$.

The ARC approach compares the (incremental) values where an arm can be chosen by sampling through the probability distribution \eqref{eq:nu} or choosing an arm with the maximum (incremental) values. (See Algorithm \ref{Alg: ARC} and \ref{Alg: ARC Index}).

\section{Batched ARC via Large Sample Theory}
\label{sec:batched setup}
 One of the main challenges in implementing the ARC approach is to ensure that the posterior distribution can be described by \eqref{eq: transition map}.  In particular, examples in \cite{ARC} only allow the observations under different arms to be correlated Gaussian or uncorrelated. Therefore, it is not clear how one can apply the ARC algorithm  to a learning problem where the observations are not Gaussian and the arms are correlated, as in the pricing example discussed in Section \ref{sec:dynamic pricing intro}. 

In this section, we will enhance the ARC algorithm to learning problems where observations $Y^{(i)}$ arrive in batches from an exponential family,  where the arms are correlated through an (unknown) parameter, in the same manner as in the generalised linear model. We assume that the observations take values in one dimension, for notational simplicity; the extension to multiple dimensions is straightforward. 

\subsection{Batched bandit problem}
Let $\{x_1, ..., x_K\} \subseteq \mathbb{R}^p$ be a collection of features  corresponding to choices $\{1,...,K\}$, and let $\bs{\theta}$ be a random variable taking values in $\Theta = \mathbb{R}^p$ representing an unknown parameter. After the agent chooses features $x_i$ at time $t$, the agent observes $Y^{(i)}_t = \big(N^{(i)}_t, (Q^{(i)}_{j,t})_{j=1}^{N^{(i)}_t}  \big)$ where $N^{(i)}_t$ is independent of each $Q^{(i)}_{j,t}$ and each $Q^{(i)}_{j,t}$ is independent with GLM density  
\begin{align}
\label{eq: exp dist}
g_0(q) \exp \Big(\bs{\Phi}_k q - G\big(\bs{\Phi}_k\big)\Big),
\end{align}
with $\bs{\Phi}_k = \phi\big(\bs{\theta}^\top x_k\big)$ for $\phi$ known.  Here, $g_0: \sR \to (0,\infty)$, $\phi : \sR \to \sR$ and $G: \sR \to \sR$ are known functions. For simplicity, we will assume that the distribution of the number of observations $(N^{(i)}_t)_{t=1}^\infty$ does not depend on the unknown parameter $\bs{\theta}$.

\subsection{Large sample theory of batched observations}
In order to apply the ARC approach to a batched bandit problem, we need to understand the evolution of the process $(M_t, D_t)$ describing the posterior distribution. As our observations arrive in batches, we use a large sample approximation to update via a classical normal-normal conjugate model for Bayesian inference. 

From \eqref{eq: exp dist}, for any $i,j,t$, the mean and variance of $Q^{(i)}_{j,t}$ given parameter $\bs{\theta} = \theta$ are 
\begin{equation}
    \label{eq:mean and var exp}
    \mathbb{E}_\theta \big(Q^{(i)}_{j,t}\big) = G'( \Phi_i ) \quad \text{and} \quad \text{Var}_\theta \big(Q^{(i)}_{j,t}\big) = G''(\Phi_i).
\end{equation}
where $\Phi_i = \phi(\theta^\top x_i)$.

Suppose that the link function $\phi$ is invertible and differentiable. If our model is non-degenerate, $G'$ must also be invertible. We define $\psi := (G' \circ \phi)^{-1}$. It then follows from the Central Limit Theorem and the Delta method that
\begin{align*}
\sqrt{n}\big(\Psi^{(i), n} - \theta^\top x_i\big) \xrightarrow[n \to \infty]{d} N\Big(0, 1/[G''\big(\Phi_i\big)\phi'(\theta^\top x_i)^2] \Big),
\end{align*}
 where $\Psi^{(i),n} = \psi\Big(\frac{1}{n}\sum_{i=1}^n Q^{(i)}_{j,t}\Big)$. 
 
Moreover,
by Slutsky's lemma, 
\begin{equation}
    \label{eq: Delta step}
\sqrt{n V(\Psi^{(i)}_n)} \big(\Psi^{(i), n} - \theta^\top x_i\big) \xrightarrow[\quad]{d} N(0, 1 ),
\end{equation}
where $V(\psi) := G''\big(\phi(\psi)\big) \big(\phi'(\psi)\big)^2$.

By large sample theory, the likelihood function for $\bs{\theta}$ when observing $\bs{\Psi}^{(i)}_t$ is

$$\pi(\bs{\theta}|\bs{\Psi}^{(i)}_t) \approx \text{constant} \times \exp \Big(-\frac{1}{2} N^{(i)}_t V(\bs{\Psi}^{(i)}_t) \big(\bs{\Psi}^{(i)}_t - \bs{\theta}^\top x_i\big)^2 \Big)$$
\noindent where the constant does not depend on $\bs{\theta}$. 

Standard Bayesian analysis shows that the posterior of $\bs{\theta}$,  after observing $\bs{\Psi}^{(i)}_t$, the prior $\cL(\bs{\theta} | \cG^U_{t-1}) = N(M_{t-1}, D_{t-1})$ can be approximately updated by
\begin{align*}
   &\pi(\bs{\theta} | \cG^U_{t-1}, \bs{\Psi}^{(i)}_t) \approx   \text{constant} \times 
   \exp \Big(-\frac{1}{2} (\bs{\theta}-M_{t})^\top D_{t}^{-1} (\bs{\theta}-M_{t}) \Big)
\end{align*}
where 
\begin{equation}
\label{eq: filtering update}
\begin{aligned}
M_{t} &= D_{t} \Big(D_{t-1}^{-1}M_{t-1} + N^{(i)}_t V(\bs{\Psi}^{(i)}_t) \bs{\Psi}_t^{(i)} x_i\Big), \\
D_{t} &= \Big(D_{t-1}^{-1} +  N^{(i)}_t V(\bs{\Psi}^{(i)}_t) x_i x_i^\top\Big)^{-1}. \\
\end{aligned}
\end{equation}

By rearranging the algebra above, we see that
\begin{equation}
\label{eq: filtering update simplified}
\begin{aligned}
M_{t} &= M_{t-1} + s^{(i)}_t D_{t-1} x_i Z^{(i)}_t,    \\
D_{t} &= D_{t-1}  - \big(s^{(i)}_t\big)^2 D_{t-1} \; x_i x_i^\top D_{t-1},
\end{aligned}
\end{equation}
where {\small $s^{(i)}_t := \sqrt{ N^{(i)}_t V(\bs{\Psi}^{(i)}_t) /(1 + N^{(i)}_t V(\bs{\Psi}^{(i)}_t) x_i^\top D_{t-1} x_i) } $} and $Z^{(i)}_t := s^{(i)}_t \big(\bs{\Psi}^{(i)}_t - M_{t-1}^\top x_i \big) \approx N(0,1)$ which can be shown using \eqref{eq: Delta step} and the predictive distribution of $\bs{\Psi}^{(i)}_t$ as in \eqref{eq: Gaussian mean propagation}.

It is worth pointing out that $(1/N^{(i)}_t) \sum_{j=1}^{N^{(i)}_t} Q^{(i)}_{j,t}$ forms a sufficient statistic for  $\bs{\theta}$ given $\big(N^{(i)}_t, (Q^{(i)}_{j,t})_{j=1}^{N^{(i)}_t}  \big)$. Therefore, as $\psi$ is invertible, the information provided by $\bs{\Psi}^{(i)}_t := \psi\Big(\frac{1}{N^{(i)}_t} \sum_{j=1}^{N^{(i)}_t} Q^{(i)}_{j,t}\Big)$ is sufficient to update the posterior.

\subsection{Overcome model degeneracy via Kalman filter}
In some settings (e.g. logistic regression), it is possible that
$(1/{N^{(i)}_t}) \sum_{j = 1}^{N^{(i)}_t} Q^{(i)}_{j,t}$ lies on the boundary (e.g. $\{0,1\}$ for logistic regression). However, the function $\psi$ is only defined on the interior (e.g.~$\psi(p) = \log \big(\frac{p}{1-p} \big)$ for logistic regression). To avoid this degeneracy, one may adopt the approach of \cite{Kalman_GLM} for Kalman filtering in the generalised linear framework to estimate $\bs{\Psi}^{(i)}_t$ using a linear expansion around $\psi^{-1}(M_{t-1}^\top x_i)$. In particular, instead of taking $\bs{\Psi}^{(i)}_t := \psi\Big(\frac{1}{N^{(i)}_t} \sum_{j=1}^{N^{(i)}_t} Q^{(i)}_{j,t}\Big)$, we may use an approximation
{\small
\begin{equation}
    \label{eq:linear interpolate}
    \bs{\Psi}^{(i)}_t \approx M_{t-1}^\top x_i + \big(\bar{Q}_t^{(i)} - \psi^{-1}(M_{t-1}^\top x_i) \big) \psi'\big( \psi^{-1}(M_{t-1}^\top x_i) \big)
\end{equation}}
where $\bar{Q}_t^{(i)} = \tfrac{1}{N^{(i)}_t} \sum_{j=1}^{N^{(i)}_t} Q^{(i)}_{j,t}$.

\subsection{ARC approach for the batched observation}
Using the approximate posterior update \eqref{eq: filtering update simplified} and \eqref{eq:linear interpolate}, we can solve the batched bandit problem using the ARC approach. The posterior dynamic \eqref{eq: transition map} can be computed using
\begin{equation}
    \label{eq:ARC posterior update}
    \begin{aligned}
    \mu_i(m,d) &= 0,  \qquad \sigma_i(m,d) = s_i(m,d) d x_i, \qquad
    b_i(m,d) &= - s_i(m,d)^2 d x_i x_i^\top d,
    \end{aligned}
\end{equation}
where $s_i(m,d) := \sqrt{{n_i V(m^\top x_i)} / ({1 + n_i V(m^\top x_i) x_i^\top d x_i})}$ and $V(\psi) := G''(\phi(\psi))\big(\phi'(\psi)\big)^2$ for $\phi$ and $G$ corresponding to the exponential family in \eqref{eq: exp dist}.

In addition to identifying the posterior dynamic through \eqref{eq: transition map} and \eqref{eq:ARC posterior update}, we need to estimate the expected reward $f(m,d)$.  Since the distribution of $Q^{(i)}_{j,t}$ depends on $\bs{\theta}$ only via $\bs{\theta}^\top x_i$, it is reasonable to assume that there exists a function $h : \sR \to \sR^K$ such that, for an observation $Y^{(i)}_t := \big(N^{(i)}_t, \big(Q^{(i)}_{j,t}\big)_{j=1}^{N^{(i)}_t} \big)$ and reward $r_i$, we have
\begin{equation}
    \label{eq:existence of h}
    \sE\big[r_i\big( Y^{(i)}_t \big) \big| \bs{\theta} \big] = h_i(\bs{\theta}^\top x_i).
\end{equation}

Using a normal approximation, we write 
\begin{equation}
\label{eq:f for batch}
\begin{aligned}
 f_i(m,d) &= \sE_{\bs{\theta} \sim N(m,d)}\big[r_i\big( Y^{(i)}_t \big) \big] \approx \int_{\sR} h_i\Big(m^\top x_i + z \sqrt{x_i^\top d x_i} \Big) \varphi(z)\d z
 \end{aligned}
\end{equation}
where $\varphi(z)$ is the standard normal density function.

  By Lemma A.4 in \cite{ARC} and \eqref{eq:ARC posterior update}, the learning premium $L$ in \eqref{eq: L expression} can be simplified to
\begin{equation}
    \label{eq:L expression_simplify}
    \begin{split}
L_i(\lambda, m,d) &= \frac{s_i(m,d)^2}{2 \lambda}  \Bigg( \sum_{j=1}^K \nu_j \big(\lambda, f(m,d) \big) g_{ij}(m,d)^2- \bigg(\sum_{j=1}^K \nu_j^\lambda \big(\lambda, f(m,d) \big) g_{ij}(m,d) \bigg)^2 \Bigg)
    \end{split}
\end{equation}
where $s_i(m,d)$ is given in \eqref{eq:ARC posterior update} and $g_{ij}(m,d) = (x_i^\top d x_j) \int_{\sR} h'_j(x_j^\top m + z ({x_j^\top d x_j})^{1/2} z)\varphi(z) \d z$. 
% For simplicity in calculation, we may estimate $g_i(m,d) \approx h_i'(x_i^\top m)$ as $d$ is small. This approximation introduces an order with negligible order compared to the approximation error considered in \cite{ARC}.
\subsection{ARC algorithm for batched bandits}
We can now describe the batched ARC algorithm.

\begin{algorithm}[H]
\label{Alg: ARC}
\SetAlgoLined
%\KwResult{Write here the result }
 \textbf{Input} $m_0,d_0, \rho, \beta$ \;
 Set $(m,d) \mapsfrom (m_0, d_0)$ \;
 \For{t = 1, 2,...}{
- Evaluate \vspace{-0.3cm}{\small $$\nu := \nu \Big(\rho \|d \|, f(m,d) + \beta(1-\beta)^{-1} L\big(\rho \|d \|, m,d\big) \Big)$$} \vspace{-0.7cm} 

via \eqref{eq:nu}, \eqref{eq:f for batch} and \eqref{eq:L expression_simplify}\;
 - Choose the  $I$-th arm by sampling $I \sim \text{Random}( \{1, 2,...,K\},  \nu)$ \; 
- Observe $Y^{(I)} = \Big(N^{(I)}_t, \big(Q^{(I)}_{j,t} \big)_{j=1}^{N^{(I)}_t} \Big)$ and collect the reward $r_{I}(Y^{(I)})$ \;
 - Update $(m,d) $ via  \eqref{eq: filtering update simplified} \;
 }
 \caption{ARC Algorithm }
\end{algorithm}

As discussed in \cite{ARC}, an alternative policy to the ARC algorithm is to choose an arm which has the maximum value of $\alpha(\lambda, m, d)$. 

\begin{algorithm}[H]
\label{Alg: ARC Index}
\SetAlgoLined
%\KwResult{Write here the result }
 \textbf{Input} $m_0,d_0, \rho, \beta$ \;
 Set $(m,d) \mapsfrom (m_0, d_0)$ \;
 \For{t = 1, 2,...}{
- Evaluate $\alpha := f(m,d) + \beta(1-\beta)^{-1} L\big(\rho \|d \|, m,d\big) $ via \eqref{eq:f for batch} and \eqref{eq:L expression_simplify}\;
 - Choose the  $I$-th arm where $I = \text{argmax} \alpha$. \; 
- Observe $Y^{(I)} = \Big(N^{(I)}_t, \big(Q^{(I)}_{j,t} \big)_{j=1}^{N^{(I)}_t} \Big)$ and collect the reward $r_{I}(Y^{(I)})$ \;
 - Update $(m,d) $ via  \eqref{eq: filtering update simplified} \;
 }
 \caption{ARC Index Algorithm }
\end{algorithm}
\section{Numerical Simulation}
\label{sec:Num sim}
\subsection{Generating environments from offline data}
\label{sec:offline data generate}
As the strategy followed determines the observations available, we cannot directly use historical data to test the algorithm. However, we can use the data to construct an environment in which to run tests. We take a Bayesian view and build a simple hierarchical model (with an improper uniform prior). Effectively, this assumes that our observations come from an exchangeable copy of the world we would deploy our bandits in, with the same (unknown) realised value of  $\bs{\theta}$. Then we will use Laplace approximation, as in \cite[Chapter 5]{Tutorial_on_Thompson_Sampling} or \cite{Laplace_approximation}, to obtain a posterior sample to simulate the markets. 
\begin{Remark}
	It is worth emphasising that, when implementing the algorithm in simulations, we do not assume that our algorithms know the distribution that we use to simulate the parameter $\bs{\theta}$. Instead simulations are initialised with an almost uninformative prior.
\end{Remark} 

To construct a posterior for $\bs{\theta}$ given historical data, we assume that we have a collection of observations  $ \cQ = \big( (q^{(1)}_j)_{j=1}^{n_1}$, $(q^{(2)}_j)_{j=1}^{n_2}$, ..., $(q^{(K)}_j)_{j=1}^{n_K} \big)$ from an exponential family modeled by \eqref{eq: exp dist}.
We denote the corresponding log-likelihood function of $\bs{\theta}$ by $\ell\big(\theta; \cQ \big)$.
%\begin{align*}
% &\ell\big(\theta; \big(q_j^{(k)}\big)\big) := \sum_{k=1}^K \Big(\phi\big(\theta^\top x^{(k)} \big)\sum_{j=1}^{n_k}q^{(k)}_j\Big)\\&\qquad - \sum_{k=1}^K n_k G\big(\theta^\top x^{(k)}\big).
%\end{align*}

Let $\hat{\theta}$ be the maximum likelihood estimator. i.e. $\hat{\theta} = \argmax_{\theta} \ell\big(\theta; \cQ \big)$. Then we may approximate the log-likelihood function by
\begin{align*}
\ell\big( \theta; \cQ \big) \approx \frac{1}{2}(\theta - \hat{\theta})^\top \del_\theta^2 \ell\big(\hat{\theta}; \cQ \big) (\theta - \hat{\theta}) + \text{constant},
\end{align*}
where the constant does not depend on $\theta$. Therefore, under the uninformative (improper, uniform) prior, we can estimate the posterior of the parameter $\bs{\theta}$ by 
\begin{equation}
\label{eq:Thetaposterior}
\bs{\theta} \sim N\Big(\hat{\theta}, -\del_\theta^2 \ell\big(\hat{\theta}; \cQ\big)^{-1} \Big).
\end{equation} 
Moreover, if the parameter for the exponential family is given in its canonical form, i.e. when $\phi$ in \eqref{eq: exp dist} is the identity, then the second derivative of the log-likelihood (i.e.~the observed Fisher information) is
$-\del_\theta^2 \ell\big(\hat{\theta}; \cQ \big) =  \sum_{i=1}^K n_i G''\big(\hat{\theta}^\top x_i\big) x_i x_i^\top.$ We can use \eqref{eq:Thetaposterior} to simulate values of $\bs{\theta}$ and use them as the (hidden) parameters to test our algorithms.

\subsection{A selection of other multi-armed bandit algorithms}
\label{sec: algo}
As a benchmark, we compare the ARC algorithm with other approaches to the multi-armed bandit problem. It is worth pointing out that the theoretical guarantees of these approaches are often provided given that the reward (not the observation) satisfies some distributional assumption, which is not the case for our setting. Nevertheless, the principles of these approaches can be extended to our setting. 

% We will compare the ARC algorithms to $\epsilon$-Greedy \cite{eps_Greedy}, Explore-then-commit (ETC) \cite{Linear_Gaussian_bandit, Lattimore_book}, Thompson Sampling \cite{Thompson_original, Tutorial_on_Thompson_Sampling}, Upper Confidence Bound (UCB, Bayes-UCB, UCB-tuned) \cite{Bayes_UCB, UCB_tuned}, Knowledge Gradient (KG)  \cite{Knowledge_Gradient} and Information-Directed Sampling (IDS) \cite{IDS, IDS_HETEROSCEDASTIC_NOISE}.  The reader may refer to \cite{ARC}, Burtini et al. \cite{adaptive_learning_survey} or the indicated references for discussion of these algorithms in more detail. 

For convenience of the reader, we give a summary of the multi-armed bandit approaches considered; further discussion can be found in \cite{ARC}. 

\begin{itemize}
	\item \textbf{$\epsilon$-Greedy (GD)} \cite{eps_Greedy}: For $\epsilon \in (0,1)$, we choose an arm based on its maximal expected reward with probability $1-\epsilon$,  and choose an arm uniformly at random and with probability $\epsilon$.
	\item \textbf{Explore-then-commit (ETC)}  \cite{Linear_Gaussian_bandit, Lattimore_book}: For $\epsilon \in (0,1)$, we choose an arm uniformly at random in the first $\lfloor \epsilon T \rfloor$ trials and then choose an arm which maximises the estimated expected reward for each remaining trial.
	
	\item \textbf{Thompson Sampling (TS)}  \cite{Thompson_original,Tutorial_on_Thompson_Sampling}: We take a posterior sample $\tilde{\bs{\theta}}$ from the prior and then choose an arm to maximise the expected reward evaluated with the sample $\tilde{\bs{\theta}}$.
	
	\item \textbf{Upper Confidence Bound (UCB)}  \cite{Bayes_UCB, GLM_bandit, UCB_tuned}: This is the optimistic strategy of choosing an arm based on an upper bound on the reward. In the classical UCB algorithm \cite{UCB_tuned}, the confidence bound is constructed by considering the reward of each arm independently. In particular, at time $t$, the classical UCB chooses an arm with the maximum
	$\bar{r}^{(i)}_t + \big(2 \log T/n^{(i)}_t\big)^{1/2}$ where $\bar{r}^{(i)}_t$ is the average reward and $n^{(i)}_t$ is the number of attempt on the $i$-th arm at time $t$ and $T$ is the horizon.
	
	Since these classical methods do not incorporate correlated learning, we will also consider the Bayes-UCB algorithm \cite{Bayes_UCB} where the arm is chosen to maximise the upper quantiles at level $p = 1-1/\big({t(\log T)^c}\big)$ where $c$ is a hyper-parameter chosen by the decision maker
	\footnote{It is not clear how to extend the UCB method considered by \cite{GLM_bandit} for the generalised linear model to our more general framework. This is because their confidence bound is constructed based on the assumption that the \emph{reward} follows the generalised linear model.}.

	\item \textbf{Knowledge Gradient (KG)} \cite{Knowledge_Gradient}: This is a Bayes one-step look ahead policy where we pretend that we will only update our posterior once, immediately after the current trial.  The decision is made by optimising this objective.
	
	\item \textbf{Information-Directed Sampling (IDS)} \cite{IDS, IDS_HETEROSCEDASTIC_NOISE}: This algorithm chooses based on the information ratio between the single-step regret and the information gain. The information gain can be chosen based on the reduction of the Shannon entropy.
\end{itemize}

\subsection{Dynamic Pricing Simulation}
We now focus on the dynamic pricing problem discussed in Section \ref{sec:dynamic pricing intro}.

In this model, at the beginning of day $t$, we choose the price of a single product from the set $\{c_1, ..., c_K\} \subseteq \sR$. With chosen price $c_i$, on day $t$, we observe $N^{(i)}_t$ customers arriving at the store. Inspired by Figure \ref{fig: price logit relation}, we model the probability that each customer buys the product using a logistic model, i.e.~the relation between the demand (probability of buying) $p(c_i)$ and the price $c_i$ is given by
$$\text{logit}\big(p(c_i)\big) = \bs{\theta}_0+ \bs{\theta}_1 c_i  = ( \bs{\theta}_0, \bs{\theta}_1)^\top (1, c_i) =: \bs{\theta}^\top x_i$$ 
where $\text{logit}(p) := \log\big(\frac{p}{1-p}\big)$, (i.e.~$G(z) = \log(1 + e^z)$ and $\phi(z) = z$ for $\phi$ and $G$ in  \eqref{eq: exp dist}). 

At the end of day $t$, we observe $Y^{(i)}_t := \big(N^{(i)}_t, \big(Q^{(i)}_{j,t}\big)_{j=1}^{N^{(i)}_t} \big)$, where $Q^{(i)}_{j,t}$ indicates whether the product (with price $c_i$) was bought by the $j$th customer of the $t$th day. The reward on day $t$ is our total revenue, given by $r_i\big(Y^{(i)}_t\big) = c_i \sum_{j=1}^{N_t^{(i)}}Q^{(i)}_{j,t}$. Using \eqref{eq:mean and var exp}, we can explicitly express the function $h_i$ in \eqref{eq:existence of h} by
 by $h_i(\cdot) = \sE[ N^{(i)}_t] c_i G'(\cdot)$.

\subsubsection{Market data and simulation environment}
\label{sec: posterior hierachi}
Dub\'e and Misra  \cite{Data_dynamic_pricing}, in stage one of their experiment, randomly assigned one of ten experimental pricing cells to 7,867 different customers who reached ZipRecruiter's paywall and estimated the price-dependence subscription rate. The exact numbers of customers that were assigned to each price are not reported. Hence, we assume that there are exactly $787$ customers for each price. We then use their reported subscription rate to estimate the number of customers who subscribed given each price.%\footnote{One can check that our estimate yields the same confidence intervals of subscription rates as shown in \cite[Figure 1]{Data_dynamic_pricing}.} 

Using this data, we apply \eqref{eq:Thetaposterior} as in Section \ref{sec:offline data generate} to infer an approximate posterior distribution given this preexisting data: 

\begin{equation}\label{eq:thetaposteriorconcrete}
\begin{aligned}
&\cL ( \bs{\theta}  | \cQ ) \approx N\Bigg( 
\begin{pmatrix}
-6.4 \times 10^{-1} \\ -4.0\times 10^{-3} 
\end{pmatrix},
\begin{pmatrix}
1.9 \times 10^{-3} & -8.9 \times 10^{-6}\\ -8.9 \times 10^{-6} & 6.8\times 10^{-8} 
\end{pmatrix}
\Bigg).
\end{aligned}
\end{equation}

% \begin{figure}[h]
% 	\centering
% 	\includegraphics[trim=6cm 1cm 7cm 1cm, width=0.7\linewidth]{price_prob_relation_simulation}
% 	\caption{(a): The (hidden) average subscription rate in our simulation, with two standard deviation bands. (b): The average reward per customer in our simulation, with two standard deviation bands.}
% 	\label{Fig: market sim}
% \end{figure}

In order to compare performance, we will consider each algorithm over a period of one year (365 days) and only allow the agent to change the price at the end of each day. We assume that a common price must be shown to all customers on each day. We also assume that the chosen price does not affect the number of customers reaching the paywall, i.e.~we assume that $N^{(i)}_t \equiv N_t$.

We run $10^4$ independent simulations where for each simulation the underlying parameter $\bs{\theta}$ is sampled from \eqref{eq:thetaposteriorconcrete}. We also independently sample $(N_t)_{t=1}^{365} \sim \text{Poisson}(270)$ to represent the number of visitors on each day. We apply the ARC approaches and other approaches with posterior parameter $(M,D)$ updating through \eqref{eq: filtering update simplified}  and \eqref{eq:linear interpolate} starting with a moderately uninformative prior $m_0 = (0,0)$ and $d_0 = I_2$. We use Monte-Carlo simulation (as discussed in \cite{ARC}) to compute relevant terms for KG and IDS. This incurs a noticeably higher computational cost than other approaches. 

% \footnote{In \cite{Data_dynamic_pricing}, they observed that ZipRecruiter.com had roughly $8,000$ visitors per month. Hence, it is reasonable to assume that there are roughly $270$ visitors per day.}.
% The simulated subscription probability and the simulated expected revenue per customer, for each price level, are illustrated in Figure \ref{Fig: market sim}.

% We shall perform the ARC approach (Algorithm \ref{Alg: ARC}) and compare to other approaches for the multi-armed bandit problem;  $\epsilon$-Greedy \cite{eps_Greedy}, Explore-then-commit (ETC) \cite{Linear_Gaussian_bandit, Lattimore_book}, Thompson Sampling \cite{Thompson_original, Tutorial_on_Thompson_Sampling}, Upper Confidence Bound (UCB, Bayes-UCB, UCB-tuned) \cite{Bayes_UCB, UCB_tuned}, Knowledge Gradient (KG)  \cite{Knowledge_Gradient} and Information-Directed Sampling (IDS) \cite{IDS, IDS_HETEROSCEDASTIC_NOISE}.  (We extend the analogy of these methods to our setting through Bayesian point of view, if necessary. The reader may refer to \cite{ARC}, Burtini et al. \cite{adaptive_learning_survey} or the indicated references for discussion of these algorithms in more detail.) Any relevant terms in these algorithms are evaluated through posterior approximation via the process $(M,D)$, as described in \eqref{eq: filtering update simplified}, starting with a moderately uninformative prior $m_0 = (0,0)$ and $d_0 = I_2$.

To assess the performance of each algorithm, we compare the cumulative pseudo-regret of each algorithm given the true parameter $\bs{\theta}$:
\[R\big(\bs{\theta} , T, (I_t)\big) := \sum_{t=1}^T \Big(\max_{i =1, ..., K}h_i(\bs{\theta}^\top x_i) - h_{I_t}(\bs{\theta}^\top x_{I_t})\Big),\]
where $h_i(y)= 270 c_i G'\big( \theta^\top x_i\big)$, $G(y) = \log(1 + e^y)$, and  $(I_t)$ is the sequence of actions that the algorithm chooses.

\subsubsection{Simulation results}

\begin{figure}[h]
 	\centering
 	\includegraphics[clip, trim=2cm 5cm 2cm 4cm, width=1\linewidth]{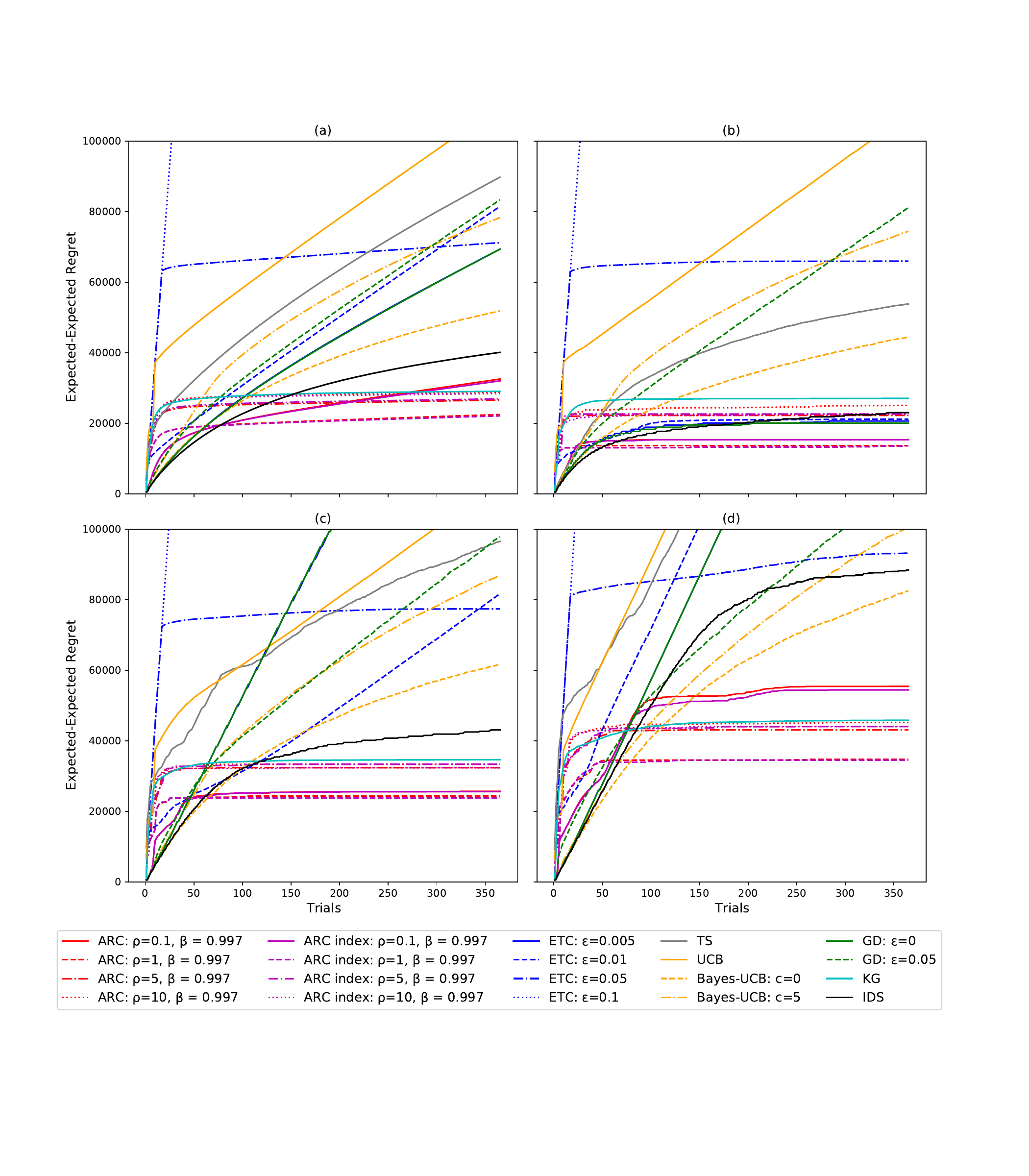}
 	\caption{(a): Mean, (b): Median, (c): 0.75 quantile and (d): 0.90 quantile of cumulative expected-expected pseudo regret}
 	\label{fig:regretplot}
 \end{figure}

Figure \ref{fig:regretplot} shows the mean, median, 0.75 quantile and 0.90 quantile of cumulative pseudo-regret (against $T$) of the ARC and other approaches. We see that most algorithms outperform the classical UCB used in \cite{bandit_and_dynamic_price}, which is unsurprising as this approach ignores correlation between demand at different prices. We also see that the ARC (and also ARC index) algorithm outperforms most algorithms both on average and in extreme cases (as shown by the quantile plots). The convergence of many well-known algorithms (e.g.~Thompson, Bayes-UCB and IDS) is too slow, as the horizon ($T = 365$) is relatively short. The performance of the ARC algorithms is a little better than the KG algorithm, but also requires lower computational cost as we do not need to use Monte-Carlo simulation to compute relevant terms. This demonstrates the benefit of studying the ARC algorithm to handle learning problems where the reward and observation are different. 

 \begin{figure}
	\centering
	\includegraphics[clip, trim=2.5cm 1cm 3cm 2cm, width=1\linewidth]{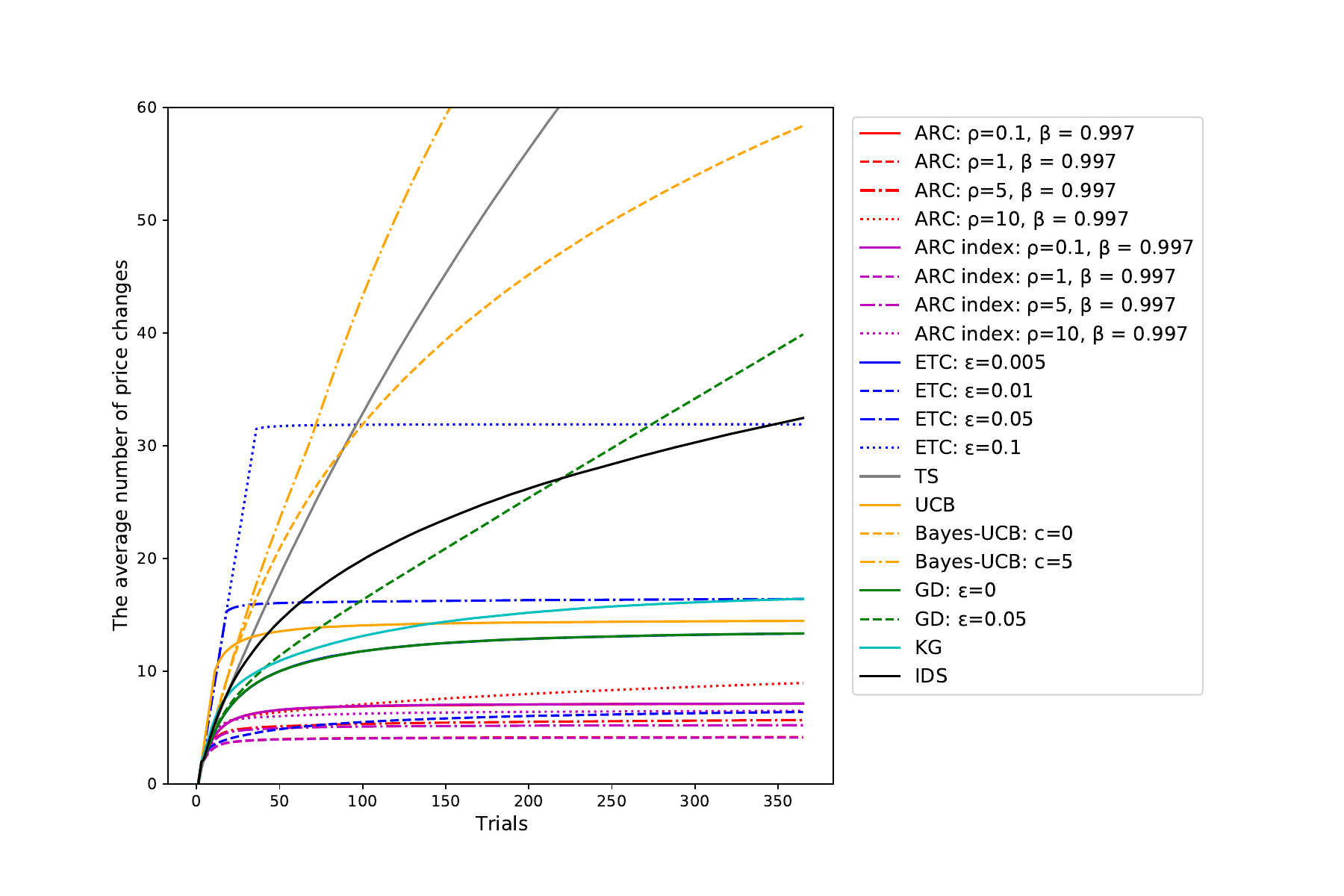}
	\caption{The average number of price changes}
	\label{Fig:Switch}
\end{figure}

In addition to the regret criteria, when displaying the price, the agent may not want to change the price too frequently. Figure \ref{Fig:Switch} shows that the ARC algorithm and KG algorithm typically require a small number of price changes, yet achieve a reasonably low regret (as shown in Figure \ref{fig:regretplot}), whereas Thompson sampling and UCB-type algorithms require a larger number of price changes.

\textbf{Acknowledgements:} Samuel Cohen also acknowledges the support of the Oxford-Man Institute for Quantitative Finance, and the UKRI Prosperity Partnership Scheme (FAIR) under the EPSRC Grant EP/V056883/1, and the Alan Turing Institute. Tanut Treetanthiploet thanks the University of Oxford for research support while completing this work, and acknowledges the support of the Development and Promotion of Science and Technology Talents Project (DPST) of the Government of Thailand and the Alan Turing Institute.

 \bibliographystyle{siam}

\bibliography{bibliography} 

\end{document}